\documentclass[10pt]{article}
\usepackage{amssymb,amsmath,amsfonts,amsthm}
    \newtheorem{theorem}{Theorem}[section]

    \newtheorem{proposition}[theorem]{Proposition}

\begin{document}
\title{On the existence of maximizers for functionals with critical
exponential growth in $\mathbb{R}^2$}
\author{Cristina Tarsi \thanks{ e-mail: tarsi@mat.unimi.it. \ } ,\\
Dipartimento di Matematica,\\ Universit\`{a} degli Studi di Milano, \\
I-20133 Milano, Italy}
\date{}
\maketitle \noindent {\footnotesize {\bf Abstract} \noindent We
investigate the problem of existence of a maximizer for
\[
S(\alpha,4\pi)=\sup_{\|u\|=1 }{\int_B\left(e^{4\pi
u^2}-1\right)|x|^{\alpha}dx,}
\]
where $B$ is the unit disk in $\mathbb{R}^2$ and $\alpha >0$. We
prove that supremum is attained for $\alpha$ small.}\vspace*{10pt}

 \medskip
%{\bf AMS subject classification:} Primary 35J60, secondary  58E05.
% main text
\section{Introduction}
Let $H^1_0(\Omega)$ be the Sobolev space over a bounded domain
$\Omega \subset \mathbb{R}^N$, with Dirichlet norm
$\|u\|^2=\int_{\Omega}|\nabla u|^2dx$. The Sobolev embedding
theorem states that $L^p(\Omega )\subset H^1_0(\Omega)$, for
$1\leq p \leq 2^{*}=\frac{2N}{N-2}$; equivalently, if we set
\[
S_{N}(p)=\sup_{\|u\|\leq 1}\int_{\Omega}|u|^pdx,
\]
then
\begin{eqnarray*}
\begin{array}{ll}
    S_{N}(p)<\infty, &  \textrm{for} \,\,\,1<p\leq 2^{*}=\frac{2N}{N-2} ; \\
    S_{N}(p)=\infty, &  \textrm{for} \,\,\,p>2^{*}; \\
\end{array}
\end{eqnarray*}
furthermore,  the value of the best Sobolev constant
$S_{N}(2^{*})$ is explicit, independent of the domain $\Omega$ and
it is known that it is never attained in any bounded smooth
domain. The maximal growth $|u|^{2^*}$ allowed is called
``critical" Sobolev growth. If $N=2$, every polynomial growth is
admitted, but it is easy to show that $H_0^1 (\Omega)\nsubseteqq
L^{\infty}(\Omega)$: in this case it is well known that the
maximal growth allowed to a function $g:\mathbb{R}\to
\mathbb{R}^+$ such that $\sup_{\|u\| \leq 1}\int g(u) <\infty$ is
of exponential type. More precisely, the Trudinger Moser
inequality states that for bounded domain $\Omega \subset
\mathbb{R}^2$
\begin{eqnarray*}
    \sup_{\|u\| \leq 1}{\int_{\Omega}e^{\gamma u^2}dx}&\leq&C(\gamma)|\Omega|
    \leq C(4\pi)|\Omega|, \,\,\,\,\textrm{for} \,\,\,\gamma\leq 4\pi ; \\
    \sup_{\|u\| \leq 1}\int_{\Omega}e^{\gamma u^2}dx&=&\infty,\,\,\,\,\textrm{for} \,\,\,\gamma>4 \pi, \\
\end{eqnarray*}
see \cite{Po}, \cite{Tr} and \cite{M}. In contrast with the
Sobolev case, the value $C(4\pi)$ is attained when $\Omega
=B_1(0)$ is the unit ball in $\mathbb{R}^2$, as proved in an
interesting paper by Carleson and Chang \cite{CC} (see also
\cite{dFdOR}). This result was extended to general bounded domains
in $\mathbb{R}^2$ by Flucher \cite{F}.

In this paper we consider the maximization problem
\begin{equation}\label{S_alpha_gamma}
S(\alpha,4\pi)=\sup_{\|u\|=1 }{\int_B\left(e^{4\pi
u^2}-1\right)|x|^{\alpha}dx,}
\end{equation}
where $\alpha >0$ and $B$ is the unit ball in $\mathbb{R}^2$. Here
we give a partial answer to a question proposed by Secchi and
Serra in a recent paper (see \cite{SS}): is the supremum attained
for any $\alpha>0$? Our main result states that $S(\alpha, 4\pi)$
 is attained, at least  if the parameter $\alpha$ is small enough.

Problem \eqref{S_alpha_gamma} can be seen as a natural
two-dimensional extension of the H\'enon-type problem
\begin{equation}\label{Henon}
\sup_{\|u\|=1
}{\left(\int_B|u|^p|x|^{\alpha}dx\right)^{2/p}}=\sup_{u\neq 0
}{\frac{\left(\int_B|u|^p|x|^{\alpha}dx\right)^{2/p}}{\int_B|\nabla
u|^2dx}}
\end{equation}
in $\mathbb{R}^N$ with $N\geq 3$ and $1<p<2^*$, which has been
widely investigated in the last few years. It is easy to verify
that \eqref{Henon} is achieved at least by a positive function;
since the quotient in \eqref{Henon} is invariant under rotations,
it is natural to ask if the supremum is achieved by a radial
function. A very interesting result obtained by Smets, Su and
Willem (\cite{SSW}) shows that a symmetry breaking phenomenon
occurs for any $p \in (2,2^*)$: in details, for every $p$ in the
subcritical range the supremum  in \eqref{Henon} is attained by a
non radial function when $\alpha \to \infty$. This result has
generated a line of research on the H\'enon-type equations (see
references in \cite{SS}). On the contrary, the H\'enon-type
problem in $\mathbb{R}^2$ with exponential nonlinearities seems to
have been much less studied. Very recently, Calanchi and Terraneo
(see \cite{CT}) proved some results about the existence of non
radial maximizers for the variational problem
\[
\sup_{\|u\|=1 }{\int_B\left(e^{\gamma u^2}-1-\gamma
u^2\right)|x|^{\alpha}dx}
\]
where $\alpha>0$ and $0<\gamma<4\pi$; in the same line is the work
by Secchi and Serra, \cite{SS}, where the authors prove a symmetry
breaking result for problem
\[
\sup_{\|u\|=1 }{\int_B\left(e^{\gamma
u^2}-1\right)|x|^{\alpha}dx.}
\]
In both papers the authors consider only subcritical exponential
growth: in this case, indeed, the existence of a maximizer can be
proved with standard arguments (see \cite{SS}, proof of
Proposition 1). On the contrary, in the critical case, that is,
when $\gamma =4\pi$, it is not clear if the supremum $S(\alpha,
4\pi)$ is attained or not. On one hand, it seems not possible to
adapt the proof suggested by Secchi and Serra, which deeply
depends on the hypothesis of subcritical growth. On the other
hand, due to the presence of the weight $|x|^{\alpha}$ in front of
the nonlinearity, problem \eqref{S_alpha_gamma} cannot be reduced
to a one dimensional problem using the technique of Schwarz
symmetrization, as proposed by Carleson and Chang.

Our result depends on a new notion of symmetrization, the so
called \emph{spherical symmetrization with respect to a measure},
which is the counterpart of Schwarz symmetrization in the
unweighted problem. Although we symmetrize with respect to a
measure $\mu$ which is different from the Lebesgue one, a result
by Schulz and Vera de Serio  (see \cite{SVS}) states that the
gradient norm does not increase, as in the classical case (the
result is valid only in $\mathbb{R}^2$ and with suitable
assumption on the measure $\mu$). This fact allows us to adapt the
proof presented by de Figueiredo, do \'O and Ruf in \cite{dFdOR},
obtaing the following result.
\begin{theorem}\label{existence_maxim}
There exists $\alpha _*>0$ such that for every $\alpha \in
(0,\alpha _*)$, $S(\alpha,4\pi)$ is attained.
\end{theorem}
We remark that the notion of symmetrization with respect to the
measure $\int_B |x|^{\alpha}dx$ gives also a geometric
interpretation of the changes of variable performed by Smets, Su
and Willem, and later by Secchi and Serra, when dealing with
radial functions: see Remark 1 at the end of Section 3.

\section{Symmetrization with respect to a measure}
In this section we recall the main definitions and properties of
symmetrization: we refer to \cite{K} or to \cite{Ba}. We start by
a review of the standard definitions. Let $\Omega$ be a bounded
domain in $\mathbb{R}^2$. We denote by $|\Omega|$ the Lebesgue
measure of $\Omega$ and by $\mathcal{L}_0(\Omega)$ the set of
Lebesgue measurable functions defined in $\Omega$ up to a.e.
equivalence. For every function $u\in \mathcal{L}_0(\Omega)$, we
define the distribution function $\phi_u$ of $u$ by the formula
\[
\phi_u(t)=|\{x\in \Omega : |u(x)|>t\}|.
\]
A measurable function $u$ in $\mathbb{R}^n$ is called
\emph{radially symmetric}, or radial, for short, if
$u(x)=\tilde{u}(r)$, $r=|x|$; it is called \emph{rearranged} if it
is nonnegative, radially symmetric and $\tilde{u}$ is a
non-increasing function of $r>0$; we also impose that
$\tilde{u}(r)$ be right-continuous. We will write $u(x)=u(r)$ by
abuse of notation. The \emph{spherical symmetric rearrangement}
$u^*$ of $u$ is the unique rearranged function defined in
$\Omega^*$ which has the same distribution function as $u$, that
is, for every $t>0$
\[
\phi_u(t)=|\{x\in \Omega : |u(x)|>t\}|=\phi_{u^*}(t)=|\{x\in
\Omega^* : |u^*(x)|>t\}|,
\]
where $\Omega^*=B_R(0)$ is the ball having the same volume as
$\Omega$, i.e. $|\Omega|=\omega_n R^n$ (here $\omega_n$ is the
volume of the unit sphere in $\mathbb{R}^n$). Then,
\begin{eqnarray}\label{sphrearr}
u^*(x)&=&\inf{\{t>0:\phi_u(t)\leq \omega_n|x|^n\}}\\
\nonumber &=& \sup{\{t>0:\phi_u(t)>\omega_n|x|^n\}}.
\end{eqnarray}
A rearranged function coincides with its spherical rearrangement.
Since the distribution functions of $u$ and $u^*$ are identical,
$\int _{\Omega}|u|^pdx=\int _{\Omega^*}{(u^*)}^pdx$ for every $p
\in [1,+\infty )$; moreover, for every nonnegative, increasing and
left-continuous real function $\Phi$
\[
\int_{\Omega}\Phi(u)dx=\int_{\Omega^*}\Phi(u^*)dx.
\]
Finally, if $u\in W^{1,p}_0(\Omega)$, then $u^*\in
W^{1,p}_0(\Omega^*)$ and
\begin{equation}\label{PS}
\int_{\Omega^*}|\nabla u^*|^pdx\leq\int_{\Omega}|\nabla u|^pdx,
\end{equation}
for $p\in (0,+\infty)$: this is the celebrated Polya-Szeg\"{o}
inequality.

As a natural generalization of the spherical symmetrization (or
Schwarz symmetrization), one can introduce the \emph{spherical
symmetrization with respect to a measure $\mu$} defined on the
domain $\Omega$. We refer to \cite{SVS}. Let
$p:\mathbb{R}^n\rightarrow \mathbb{R}^+$ be a nonnegative,
measurable and locally integrable function, and consider the
absolutely continuous measure $\mu$ given by
\[
\mu(A)=\int_A pdx
\]
for any Lebesgue measurable set $A$ in $\mathbb{R}^n$. The
distribution function $\phi_{\mu,u}$ of $u$ with respect to the
measure $\mu$ is given by
\[
\phi_{\mu,u}(t)=\mu\left(\{x\in \Omega : |u(x)|>t\}\right);
\]
as in the classical case, $\phi_{\mu,u}$ is a monotone,
non-increasing and right continuous function. The \emph{spherical
symmetric rearrangement} $u^*_{\mu}$ of $u$ \emph{with respect to
the measure} $\mu$ is the unique rearranged function defined in
$\Omega^*_{\mu}$ whose (classical) distribution function is the
same as the distribution function (with respect to the measure
$\mu$) of $u$; that is, for every $t>0$
\[
\phi_{\mu,u}(t)=\mu\left(\{x\in \Omega :
|u(x)|>t\}\right)=\phi_{u^*}(t)=|\{x\in \Omega^*_{\mu} :
|u^*(x)|>t\}|,
\]
where $\Omega^*_{\mu}=B_R(0)$ is the ball centered at the origin
with $\mu(\Omega)=|\Omega^*_{\mu}|=\omega_n R^n$. Then,
\begin{eqnarray}\label{sphrearrmu}
u^*_{\mu}(x)&=&\inf{\{t>0:\phi_{\mu,u}(t)\leq \omega_n|x|^n\}}\\
\nonumber &=& \sup{\{t>0:\phi_{\mu,u}(t)>\omega_n|x|^n\}}.
\end{eqnarray}
Obviously, the spherical symmetric rearrangement
$u^*_{\mathcal{L}}$ with respect to the Lebesgue measure
 is the classical symmetric rearrangement by Schwarz.
However, if $\mu$ is not the Lebesgue measure, a rearranged
function, in the sense defined above, will not coincide with its
$\mu$-rearrangement $u^*_{\mu}$, since an extra
contraction/dilation will take place. In particular, if we
consider the density function
$p_{\alpha}(x)=|x|^{\alpha}:\mathbb{R}^n\rightarrow \mathbb{R}^+$
with $\alpha>0$ and the associated measure
\begin{equation}\label{mualpha}
\mu_{\alpha}(A)=\int_A|x|^{\alpha}dx
\end{equation}
defined on the unit sphere $B$ in $\mathbb{R}^n$, then the
$\mu_{\alpha}$-rearrangement of a rearranged function $u(r)$ (that
is, of a non-negative, radial and non-increasing function $u$) is
defined by the formula
\begin{equation}\label{radialmurearr}
u^*_{\alpha}(r)=u\left(r^{\frac n{\alpha+n}}\sqrt[\alpha
+n]{\frac{\alpha +n}n}\right),
\end{equation}
where $r\in B\left(0,\sqrt[n]{\frac{\alpha +n}n}\right)$.

As in the classical case, for every nonnegative, increasing and
left-continuous real function $\Phi$
\[
\int_{\Omega}\Phi(u)d{\mu}=\int_{\Omega^*_{\mu}}\Phi(u^*_{\mu})dx,
\]
so that
$\|u\|_{L^p(\Omega,\mu)}=\|u^*_{\mu}\|_{L^p(\Omega^*_{\mu},\mathcal{L})}$
for every $p\in [1,+\infty)$. As regard the gradient norm, it is
not known (to our knowledge) if the Polya-Szeg\"{o} inequality can
be maintained for all $\mu$-rearrangements and for $p>0$, $n\geq
1$. With the assumptions stated above on $\mu$, if $u\in
W^{1,1}(\mathbb{R}^n)$, then $u^*_{\mu}\in W^{1,1}(\mathbb{R}^n)$;
furthermore, Schulz and Vera de Serio have proved the following
result:
\begin{theorem}[F. Schulz, V. Vera de Serio]\label{SVS}
Let $p\in\mathcal{C}^0(\bar{D})$ be a nonnegative function on a
simply-connected domain $D$ such that $\log{p}$ is subharmonic
where $p>0$; suppose that $u\in W^{1,2}(\mathbb{R}^2)$ is a
non-negative function with compact support in $D$. Then $u^*_{\mu}
\in W^{1,2}(\mathbb{R}^2)$, and the inequality
\begin{equation}\label{PSmu}
\int_{\mathbb{R}^2}|\nabla
u^*_{\mu}|^2dx\leq\int_{\mathbb{R}^2}|\nabla u|^2dx
\end{equation}
holds.
\end{theorem}
We remark that Theorem \ref{SVS} states that the gradient of the
$\mu$-rearrangement does not increase in the $L^2(\mathbb{R}^2)$
norm (that is, considering $\mathbb{R}^2$ endowed with the
Lebesgue measure); different results can be found in \cite{Ta} and
in \cite{SW} where a similar inequality is obtained for the
$L^2(\mathbb{R}^2, \mu)$ norms.
\section{Existence of a maximizer for $S(\alpha, 4\pi)$}
This section is devoted to the proof of our main result, Theorem
\ref{existence_maxim}. As well known, in the ``unweighted" case
$\alpha=0$ the supremum $S(0,4\pi)$ is attained: this is the
celebrated result due to Carleson and Chang \cite{CC}. In the
subcritical case $\gamma <4\pi$, $S(\alpha,\gamma)$ is still
attained, as pointed out by Serra and Secchi in \cite{SS} (and the
proof is quite easy), whereas in the supercritical case $\gamma >
4\pi$, $S(\alpha,\gamma)=+\infty$ for every $\alpha>0$, as proved
by Calanchi and Terraneo \cite{CT} testing with a suitable
sequence of (radial) functions.

The critical case $\gamma = 4\pi$ is more delicate. If we consider
the radial version of the maximization problem
\eqref{S_alpha_gamma}, that is,
\begin{equation}\label{S_alpha_gamma_rad}
S^{\textrm{rad}}(\alpha,\gamma)=\sup_{\begin{array}{c}
  u\in H_{0,\textrm{rad}}^1(B) \\
  \|u\|= 1 \\
\end{array}}{\int_B\left(e^{\gamma u^2}-1\right)|x|^{\alpha}dx,}
\end{equation}
it is not hard to prove that the problem is still ``subcritical",
provided that $\gamma<4\pi +2\pi \alpha$, as proved by Secchi and
Serra in \cite{SS}. More in details, they proved that
\[
S^{\textrm{rad}}(\alpha, 4\pi)=\frac2{\alpha +2}S\left(0,4\pi
\frac2{\alpha +2}\right)
\]
and standard arguments show that $S(0,4\pi \frac2{\alpha +2})$ is
actually attained by a radial function. See also the remark at the
end of Section 3 in \cite{CT}.

On the contrary, it seems not possible to reduce the problem of
maximization of $S(\alpha, 4\pi)$ in the general case, that is,
considering also non radial functions, to a subcritical one. Our
proof follows the same idea of the one given by de Figueiredo, do
\'{O} and Ruf in \cite{dFdOR} (which differs from the original
proof of Carleson and Chang by the use of the
concentration-compactness principle). Here is a short outline of
the proof:
\begin{itemize}
    \item if $S(\alpha, 4\pi)$ is not attained, then by the
    concentration-compactness alternative of P.L. Lions there is a
    normalized maximizing and concentrating sequence $v_n$;
    \item by means of symmetrization with respect to the measure
    $\mu_{\alpha}=\int |x|^{\alpha} dx$, one can prove an upper
    bound for any normalized concentrating sequences $u_n$:
    \[
     \overline{\lim}_{n\rightarrow
+\infty}{\int_{B^*_{\alpha}}\left(e^{4\pi |u_n|^2}-1
\right)dx}\leq \frac2{\alpha +2}\pi e
    \]
    \item give an explicit function $\omega \in H^1_0(B)$ such
    that $\|\omega\|=1$ and
    \[
\int_B\left(e^{4\pi \omega^2}-1\right)|x|^{\alpha} dx >
\frac2{\alpha +2}\pi e.
\]
\end{itemize} It is clear, then, that the notion of spherical symmetrization
 with respect to a measure is the fundamental tool which allows
 to reduce the weighted problem $S(\alpha, 4\pi)$ to a one dimensional problem.
 See also the remarks at the end of the proof.

First of all, let us recall the concentration-compactness result
by P.L. Lions \cite{L} (adapted to the 2-dimensional case):
\begin{proposition}[P.L. Lions]\label{conc_comp}
Let $\Omega$ be a bounded domain in $\mathbb{R}^2$, and let
$\{u_n\}$ be a sequence in $H_0^1(\Omega)$ such that
$\|u_n\|_{H_0^1}\leq 1$ for all $n$. We may suppose that $u_n
\rightharpoonup u$ weakly in $H_0^1(\Omega)$, $|\nabla
u_n|^2\rightarrow \nu$ weakly in measure. Then either
\newline(i) $\nu= \delta _{x_0}$, the Dirac measure of mass $1$
concentrated at some $x_0 \in \bar{\Omega}$, and $u\equiv
0$,\newline or \newline (ii) there exists $\beta > 4\pi$ such that
the family $v_n=e^{u_n^2}$ is uniformly bounded in
$L^{\beta}(\Omega)$ and thus $\int_{\Omega}e^{4\pi u_n^2}
\rightarrow \int_{\Omega}e^{4\pi u^2}$ as $n\rightarrow +\infty$.
In particular, this is the case if $u$ is different from $0$.
\end{proposition}
\vspace{10pt} \textbf{Proof of Theorem \ref{existence_maxim}}
\begin{proof}
We follow \cite{dFdOR}. We say that a sequence $\{u_n\} \subset
H_0^1(B)$ is a \emph{normalized concentrating sequence} if
\newline $i)$ $\|u_n\|_{H_0^1}= 1$
\newline $ii)$ $u_n \rightharpoonup 0$ weakly in $H_0^1(B)$
\newline $iii)$ $\exists \; x_0\in B$ such that $\forall \rho>0$,  %
$\int_{B\setminus B_{\rho}(x_0)}|\nabla u_n|^2dx\rightarrow
0$.\vspace{10pt}

Let us suppose that $\{u_n\}$, $\|u_n\|= 1$, is a maximizing
sequence for \eqref{S_alpha_gamma}, that is, $\lim_{n\rightarrow
+\infty}{\int_B(e^{4\pi u_n^2}-1)|x|^{\alpha}dx}=S(\alpha, 4\pi)$.
Then, by the concentration-compactness alternative of P.L. Lions,
either $\{u_n\}$ is a normalized concentrating (and maximizing)
sequence, or $S(\alpha, 4\pi)$ is attained. To conclude the proof,
we proceed by the following steps: \vspace{5pt}\newline 1) if
$\{u_n\}$ is any normalized concentrating sequence in $H_0^1(B)$,
then
\begin{equation}\label{upperbound_alpha}
\lim_{n\rightarrow +\infty}{\int_B \left(e^{4\pi
u_n^2}-1\right)|x|^{\alpha}dx}\leq \frac2{\alpha +2}\pi e;
\end{equation}
\vspace{5pt}\newline 2) give an explicit function $\omega \in
H_0^1(B)$ such that
\[
\int_B\left(e^{4\pi \omega^2}-1\right)|x|^{\alpha}dx >
\frac2{\alpha +2}\pi e.
\]
\textbf{1) Upper bound.} Using the notion of spherical
symmetrization with respect to the measure $\mu
_{\alpha}=\int_B|x|^{\alpha}$ introduced in Section 2, and Theorem
\ref{SVS} of Schulz-Vera de Serio, it suffices to show that
\[
\overline{\lim}_{n\rightarrow
+\infty}{\int_{B^*_{\alpha}}\left(e^{4\pi |u^*_{\alpha,n}|^2}-1
\right)dx}\leq \frac2{\alpha +2}\pi e
\]
where $\{u^*_{\alpha,n}\}$ is the rearranged sequence of ${u_n}$,
with $\|u^*_{\alpha,n}\| \leq \|u_n\| = 1$, and $B^*_{\alpha}$ is
the ball centered in 0 such that $|B^*_{\alpha}|=\mu_{\alpha}(B)$,
that is,
\[
B^*_{\alpha}=B\left(0, \sqrt{\frac2{\alpha +2}}\right).
\]
Let us set $z_n=\frac{u^*_{\alpha,n}}{\|u^*_{\alpha,n}\|}$; then
\[
\int_{B^*_{\alpha}}\left(e^{4\pi |u^*_{\alpha,n}|^2}-1
\right)dx\leq \int_{B^*_{\alpha}}\left(e^{4\pi z^2_n}-1 \right)dx,
\]
so that is suffices to prove that for any radial normalized
concentrating sequence in $B(0, \sqrt{\frac2{\alpha +2}})$ the
upper bound \eqref{upperbound_alpha} holds. First, we perform a
change of variable to reduce the domain to the unit ball. Let $R
=\sqrt{\frac2{\alpha +2}} \rho$, and
$y_n(\rho)=z_n(\sqrt{\frac2{\alpha +2}} \rho)$; then
\[
2\pi\int_0^{\sqrt{\frac2{\alpha +2}}}\left(e^{4\pi z^2_n}-1
\right)RdR=2\pi \frac2{\alpha +2} \int_0^1\left(e^{4\pi
 y^2_n}-1 \right)\rho d\rho
\]
and
\[
1=\int_{B^*_{\alpha}}|\nabla z_n^*|^2dx=2\pi \int_0^1 |y_n'|^2
\rho d\rho.
\]
The proof now reads exactly as in \cite{dFdOR} (proof of Theorem
4, step 1), so we can omit it. See also \cite{CC}.
\vspace{10pt}\newline \textbf{ 2) An explicit function.} In this
step we exhibit an explicit function $\omega(x)$ such that
\[
\int_B\left(e^{4\pi \omega^2}-1\right)|x|^{\alpha}dx >
\frac2{\alpha +2}\pi e;
\]
since, by step 1), any maximizing sequence (if exists) must
satisfy $S(\alpha, 4\pi)=\lim_{n\rightarrow
+\infty}{\int_B(e^{4\pi u_n^2}-1)|x|^{\alpha}dx}\leq \frac2{\alpha
+2}\pi e$, we can conclude that $S(\alpha, 4\pi)$ is attained.
From now on we assume that $u$ is a generic radial function, and
set
\begin{equation}\label{epsilon}
\varepsilon=\frac2{\alpha +2}.
\end{equation}
As in \cite{SS}, following an idea of Smets, Su and Willem, define
the new function
\begin{equation}\label{v}
v(\rho)=\frac1{\sqrt{\varepsilon}}u(\rho ^{\varepsilon});
\end{equation}
then
\begin{equation}\label{gradv}
\int_B|\nabla u|^2dx=2\pi \int_0^1|u'|^2rdr=2\pi
\int_0^1|v'|^2\rho d\rho
\end{equation}
and,
\begin{equation}\label{e^4pi_eps_v}
\int_B\left(e^{4\pi u^2-1}\right)|x|^{\alpha}dx=2\pi \varepsilon
\int_0^1\left(e^{4\pi \varepsilon v^2-1}\right) \rho d\rho.
\end{equation}
We can now perform the change of variable introduced by Moser
\cite{M}, which transform the radial integral on $[0,1)$ into an
integral on the half-line $[0, +\infty)$,
\[
\rho=
e^{-t/2}\;\;\;\;\;\;\textrm{and}\;\;\;\;\;w(t)=\sqrt{4\pi}v(\rho);
\]
we obtain (recalling the definition \eqref{epsilon})
\begin{equation}\label{int_explicit}
\int_B\left(e^{4\pi u^2}-1\right)|x|^{\alpha}dx=\pi \frac2{\alpha
+2} \left(\int_0^{+\infty}  e^{\frac2{\alpha +2}w^2-t}dt -1\right)
\end{equation}
with
\[
\int_B|\nabla u|^2dx=\int_0^{+\infty}|w'(t)|^2dt.
\]
Following \cite{CC}, take $w:[0,+\infty) \rightarrow \mathbb{R}$
to be
\[
w(t)=\left\{ \begin{array}{ll} \frac12 t & \textrm{if}\;\; 0\leq
t\leq 2\\
\sqrt{t-1} &\textrm{if}\;\; 2\leq
t\leq 1+e^2\\
e &\textrm{if}\;\;  t\geq 1+e^2
 \end{array} \right. .
\]
Then, by direct inspection
\[
\int_0^{+\infty}|w'(t)|^2dt=1
\]
and
\begin{eqnarray}\label{stima_int_expl}\nonumber
\int_0^{+\infty}  e^{\frac2{\alpha +2}w^2-t}dt &=& \int_0^2
e^{\frac2{\alpha +2}\frac{t^2}4-t}dt%
+e^{-\frac2{\alpha +2}}\int_2^{1+e^2}e^{-(1-\frac2{\alpha +2})t}dt
\\ \nonumber && \;\;\;\; +e^{e^2 \frac2{\alpha +2}}\int_{1+e^2}^{+\infty}e^{-t}dt
\\ \nonumber  &=& \int_0^2
e^{\frac2{\alpha +2}\frac{t^2}4-t}dt \\
\nonumber && + \frac1{e}\left[ e^{-\frac{\alpha}{\alpha +2}e^2} -
\frac{\alpha +2}{\alpha} e^{-\frac{\alpha}{\alpha +2}e^2} +
\frac{\alpha
+2}{\alpha}e^{-\frac{\alpha}{\alpha +2}} \right]\\
 &=& \int_0^2 e^{\frac2{\alpha +2}\frac{t^2}4-t}dt
+A_{\alpha}
\end{eqnarray}
where
\[
A_{\alpha}=\frac1{e}\left[ - \frac2{\alpha}
e^{-\frac{\alpha}{\alpha +2}e^2} + \frac{\alpha
+2}{\alpha}e^{-\frac{\alpha}{\alpha +2}} \right]
\]
Let us now estimate the right hand side of
\eqref{stima_int_expl}, when $\alpha \rightarrow 0$. Set
$s=-\frac{t}{\alpha +2}+1$ in the integral term; then
\begin{eqnarray}
\nonumber \int_0^2 e^{\frac2{\alpha +2}\frac{t^2}4-t}dt&=& (\alpha
+2) e^{-\frac{\alpha +2}2} \int_{\frac{\alpha}{\alpha +2}}^1
e^{\frac{\alpha +2}2 s^2}ds \\
\nonumber &>& (\alpha +2) e^{-\frac{\alpha +2}2}
\int_{\frac{\alpha}{\alpha +2}}^1 e^{s^2}(1+\frac{\alpha}2s^2)ds\\
\nonumber &=& (\alpha +2) e^{-\frac{\alpha +2}2}(1-\frac{\alpha}4)
\int_{\frac{\alpha}{\alpha +2}}^1 e^{s^2}ds + \\
\nonumber && (\alpha +2)e^{-\frac{\alpha +2}2}\frac{\alpha}4
\left\{ e-\frac{\alpha}{\alpha
+2}e^{(\frac{\alpha}{\alpha+2})^2}\right\}
\end{eqnarray}
\begin{eqnarray}\label{I_1}
&=& (\alpha +2) e^{-\frac{\alpha +2}2}(1-\frac{\alpha}4) \int_0^1
e^{s^2}ds \\
\nonumber&& - (\alpha +2) e^{-\frac{\alpha +2}2}(1-\frac{\alpha}4)
\int_0^{\frac{\alpha}{\alpha +2}} e^{s^2}ds + B_{\alpha},
\end{eqnarray}
where
\[
B_{\alpha}=(\alpha +2)e^{-\frac{\alpha +2}2}\frac{\alpha}4 \left\{
e-\frac{\alpha}{\alpha +2}e^{(\frac{\alpha}{\alpha+2})^2}\right\}
\]
When $\alpha \rightarrow 0$,
\[
(\alpha +2) e^{-\frac{\alpha +2}2}(1-\frac{\alpha}4)=\frac2{e}
+\textrm{o}(1),
\]
\[
\int_0^{\frac{\alpha}{\alpha +2}} e^{s^2}ds<\frac{\alpha}{\alpha
+2}e^{(\frac{\alpha}{\alpha +2})^2}=\textrm{o}(1)
\]
and
\[
B_{\alpha}=\frac{\alpha}2 + \textrm{o}(\alpha).
\]
Combining \eqref{I_1} with the last estimates yields
\begin{equation}\label{stimaI_1}
\int_0^2  e^{\frac2{\alpha +2}w^2-t}dt=\frac2{e}\int_0^1
e^{s^2}ds+\textrm{o}(1).
\end{equation}
On the other hand,
\begin{eqnarray}\label{A_alpha}\nonumber
A_{\alpha}&=&\frac1{e}\left[ - \frac2{\alpha}
e^{-\frac{\alpha}{\alpha +2}e^2} + \frac{\alpha
+2}{\alpha}e^{-\frac{\alpha}{\alpha +2}} \right]\\
&=& e +\textrm{o}(1)\,\,\,\textrm{as}\,\,\alpha \to 0
\end{eqnarray}
Inserting \eqref{stimaI_1} and \eqref{A_alpha} in
\eqref{stima_int_expl} we obtain
\[
\int_0^{+\infty}  e^{\frac2{\alpha +2}w^2-t}dt = e
+\frac2{e}\int_0^1 e^{s^2}ds +\textrm{o}(1);
\]
therefore,  if $\omega$ is the radial function which corresponds
to $w(t)$, by \eqref{int_explicit} we have
\begin{eqnarray}\nonumber
\int_B\left(e^{4\pi \omega ^2}-1\right)|x|^{\alpha}dx&=&\pi
\frac2{\alpha
+2}\left(e +\frac2{e}\int_0^1 e^{s^2}ds -1 +\textrm{o}(1)\right)\\
\nonumber&>&  \frac2{\alpha +2} \pi e\hspace{20 pt} \textrm{when
}\alpha \rightarrow 0,
\end{eqnarray}
since $\frac2{e}\int_0^1 e^{s^2}ds>1$, as one can verify
estimating the integral with lower Riemann sum, as in \cite{CC}
(the value obtained is $\frac2{e}\int_0^1 e^{s^2}ds\approx
\frac{2.723}{e}>1$), or expanding the integrand in power series,
as in \cite{SS} (here $\frac2{e}\int_0^1 e^{s^2}ds\approx
\frac{2.906}{e}>1$).
\end{proof}
\textbf{Remark 1.} The notion of symmetrization with respect to
the measure $\mu_{\alpha}=\int|x|^{\alpha}$ is a fundamental tool
in the proof of Theorem \ref{existence_maxim}, as remarked in the
introduction, since it allows to reduce the variational problem to
a one-dimensional problem, as in the unweighted case $\alpha =0$.
Furthermore, it gives a geometric interpretation of the change of
variable \eqref{v}, originally introduced by Smets, Su and Willem
in \cite{SSW}, which allows to reduce the weighted integral
$\int_B(e^{4\pi u^2}-1)|x|^{\alpha}dx$ to the unweighted integral
$\varepsilon \int_B(e^{4\pi \varepsilon u^2}-1)dx$ if $u$ is a
radial function. Indeed, let us consider a rearranged function
$u(r)$; then, by \eqref{radialmurearr} of the previous section
(and using the notation \eqref{epsilon} for simplicity)
\[
u^*_{\alpha}(r)=u\left(r^{\varepsilon}\varepsilon
^{-\frac{\varepsilon}2}\right),
\]
so that
\begin{eqnarray}\nonumber
\int_B|\nabla u|^2dx&=&2\pi \int_0^1 |u'(s)|^2s ds\\
\nonumber &=& \frac{2\pi }{\varepsilon}
\int_0^{\sqrt{\varepsilon}}|u^{*'}_{\alpha}(r)|^2 r dr
\end{eqnarray}
and
\begin{eqnarray}\nonumber
\int_B \left(e^{4\pi u^2}-1\right)|x|^{\alpha}dx &=& 2\pi \int_0^1
\left(e^{4\pi u^2}-1\right)s^{\alpha +1}ds\\
\nonumber &=& 2\pi \int_0^{\sqrt{\varepsilon}} \left(e^{4\pi
|u^*_{\alpha}|^2}-1\right) rdr.
\end{eqnarray}
Now, set $r= \sqrt{\varepsilon}\rho$ and
$v(\rho)=u^*_{\alpha}(\sqrt{\varepsilon}\rho)$; then
\begin{equation}\label{SS1}
\frac{2\pi}{\varepsilon}
\int_0^{\sqrt{\varepsilon}}|u^{*'}_{\alpha}(r)|^2 r dr=
\frac{2\pi}{\varepsilon}\int_0^1|v'(\rho)|^2 \rho d\rho
\end{equation}
and
\begin{equation}\label{SS2}
 2\pi \int_0^{\sqrt{\varepsilon}}
\left(e^{4\pi |u^*_{\alpha}|^2}-1\right) rdr=2\pi \varepsilon
\int_0^1 \left(e^{4\pi v^2}-1\right) \rho d\rho.
\end{equation}
Equalities \eqref{SS1} and \eqref{SS2} can be restated as
\begin{eqnarray}\nonumber
\int_B|\nabla u|^2dx&=& 2\pi
\frac1{\varepsilon} \int_0^1|v'|^2 \rho d\rho,\\
\nonumber \int_B \left(e^{4\pi u^2}-1\right)|x|^{\alpha}dx&=&2\pi
\varepsilon \int_0^1 \left(e^{4\pi v^2}-1\right) \rho d\rho,
\end{eqnarray}
where
\[
v(\rho)=u^*_{\alpha}(\sqrt{\varepsilon}
\rho)=u(\rho^{\varepsilon});
\]
this is exactly the change of variable introduced by Smets, Su and
Willem in \cite{SSW}, and differs from \eqref{v} by a dilation
factor. Therefore, the change of variable $\rho = r^{\varepsilon}$
performed to obtain asymptotic estimates for the radial supremum
$S^{\textrm{rad}}(\alpha, 4\pi)$ in \cite{SS} (respectively
$S^{\textrm{rad}}(\alpha, p)$ in \cite{SSW}) coincides with the
spherical symmetrization with respect to the measure
$\mu_{\alpha}$ (rescaled so to reduce the symmetrized domain
$B^*_{\alpha}$ to the unit ball $B$).

\textbf{Remark 2.} Note that we have proved step 2 testing with a
radial function. It easy to show that if $\alpha \rightarrow
+\infty$, a function $w(x)$ such that
\[
\int_B  \left(e^{4\pi w^2}-1\right)|x|^{\alpha}dx > \frac2{\alpha
+2} \pi e,
\]
if exists, must be \emph{non radial}. Indeed, for any $u$ radial
function in $H^{1, \textrm{rad}}_0(B)$, $\int_B \left(e^{4\pi
u^2}-1\right)|x|^{\alpha}dx=2\pi \varepsilon \int_0^1 (e^{4\pi
\varepsilon v^2}-1)rdr$ by \eqref{e^4pi_eps_v}; but
\begin{eqnarray}\label{partial_int_alpha}\nonumber
\frac{\partial}{\partial \varepsilon}\left(2\pi \int_0^1 (e^{4\pi
\varepsilon v^2}-1)rdr\right)&=&\frac1{\varepsilon}2\pi \int_0^1
4\pi \varepsilon v^2e^{4\pi \varepsilon v^2}rdr
\\ &>& \frac1{\varepsilon}2\pi \int_0^1 \left(e^{4\pi
\varepsilon v^2}-1\right)rdr,
\end{eqnarray}
since $te^t>e^t-1$ for every $t>0$ (note that the inequality is
strict, and there is equality if and only if $t=0$). Integrating
the previous inequality yields
\begin{eqnarray*}
2\pi \int_0^1 (e^{4\pi v^2}-1)rdr &>& \frac1{\varepsilon }2\pi
 \int_0^1 (e^{4\pi \varepsilon v^2}-1)rdr \\
 &=&\frac1{\varepsilon ^2} \int_B (e^{4\pi u^2}-1)|x|^{\alpha}dx.
\end{eqnarray*}
Therefore, if there exists a radial function $w(x)$, with
$\|w\|\leq 1$, such that $\int_B  (e^{4\pi w^2}-1)|x|^{\alpha}dx >
\frac2{\alpha +2} \pi e= \varepsilon \pi e$, we have also
\[
S(0, 4\pi)\geq  \int_B (e^{4\pi w^2}-1)dx> \frac1{\varepsilon ^2}
\int_B (e^{4\pi w^2}-1)|x|^{\alpha}dx>\frac1{\varepsilon} \pi e;
\]
this implies that $S(0, 4\pi)$ is unbounded  as $\alpha \to
+\infty$, that is a contradiction. Hence, the problem of the
existence of a maximizer for $S(\alpha, 4\pi)$, when
$\alpha>\alpha_*$, is reduced to finding a non-radial function $w$
satisfying the lower bound $\int_B  (e^{4\pi w^2}-1)|x|^{\alpha}dx
> \frac2{\alpha +2} \pi e$ (trying to adapt the proof presented
here).

\vspace{10pt} \textbf{Remark 3.} It remains an open problem
whether the supremum $S(\alpha, 4\pi)$ is attained \emph{for
every} $\alpha>0$.

\end{document}